\documentclass{gtart}

\input gtspec
\volumenumber{3}\papernumber{5}\volumeyear{1999}
\pagenumbers{119}{135}\published{29 May 1999}
\proposed{Robion Kirby}\seconded{Walter Neumann, Yasha Eliashberg}
\received{23 November 1998}
\accepted{20 May 1999}

\usepackage{amsmath,amsfonts}
\usepackage{pstricks,pst-node}
\let\Bbb\mathbb

\newtheorem{theorem}{Theorem}
\newtheorem{proposition}{Proposition}
\newtheorem{conjecture}{Conjecture}
\newtheorem{lemma}{Lemma}

\newcommand{\ie}{that is}
\renewcommand{\d}{\partial}

\newcommand{\R}{{\Bbb R}}

\newcommand{\vf}{\vec{f}}
\newcommand{\ve}{\vec{e}}
\newcommand{\vj}{\vec{\jmath}}
\newcommand{\vv}{\vec{v}}

\newcommand{\vx}{\vec{x}}
\newcommand{\vy}{\vec{y}}

\newcommand{\heart}{\heartsuit}
\renewcommand{\diamond}{\diamondsuit}
\newcommand{\dual}{\circ}
\newcommand{\st}{{\bigm|}}

\newcommand{\Tr}{\mbox{Tr}}

\newcommand{\Wedge}{{\textstyle{\bigwedge}}}

\DeclareMathOperator{\GL}{GL}
\DeclareMathOperator{\Hom}{Hom}
\DeclareMathOperator{\so}{so}
\DeclareMathOperator{\SO}{SO}
\DeclareMathOperator{\Vol}{Vol}
\DeclareMathOperator{\vecvol}{\stackrel{\longrightarrow}{Vol}}
\newcommand{\tensor}{\otimes}

\newcommand{\Co}[1]{Conjecture~\ref{#1}}
\newcommand{\Thm}[1]{Theorem~\ref{#1}}
\newcommand{\Prop}[1]{Proposition~\ref{#1}}
\newcommand{\Sec}[1]{Section~\ref{#1}}
\newcommand{\Fig}[1]{Figure~\ref{#1}}
\newcommand{\Lem}[1]{Lemma~\ref{#1}}

\newcommand{\bracket}[1]{{\langle #1\rangle}}

\newcounter{fignum}
\newenvironment{fullfigure}[2]
    {\begin{figure}[htb]\small\begin{center}\def\fullfiga{#1}\def\fullfigb{#2}}
    {\caption{\fullfigb}\label{\fullfiga}\end{center}\end
{figure}}

\newgray{silver}{.9}

\psset{linewidth=.4pt,dash=3pt 3pt,doublesep=.05,dotsize=1pt 5}
\SpecialCoor

\begin{document}

\title{The bottleneck conjecture}
\author{Greg Kuperberg}
\address{Department of Mathematics,
University of California\\
One Shields Avenue, Davis, CA 95616, USA}
\email{greg@math.ucdavis.edu}
\begin{abstract}
The Mahler volume of a centrally symmetric convex body $K$ is defined as $M(K)
= (\Vol K)(\Vol K^\dual)$.  Mahler conjectured that this volume is minimized
when $K$ is a cube.  We introduce the bottleneck conjecture, which stipulates
that a certain convex body $K^\diamond \subset K \times K^\dual$ has least
volume when $K$ is an ellipsoid.   If true, the bottleneck conjecture would
strengthen the best current lower bound on the Mahler volume due to Bourgain
and Milman.  We also generalize the bottleneck conjecture in the
context of indefinite orthogonal geometry and prove some special cases
of the generalization.

This article is in the \texttt{xxx} archive as:\qua
\textbf{\texttt{math.MG/9811119}}
\end{abstract}
\asciiabstract{The Mahler volume of a centrally symmetric convex body K
is defined as M(K)= (Vol K)(Vol K^dual).  Mahler conjectured that this 
volume is minimized
when K is a cube.  We introduce the bottleneck conjecture, which stipulates
that a certain convex body $K^diamond subset K X K^dual has least
volume when K is an ellipsoid.   If true, the bottleneck conjecture would
strengthen the best current lower bound on the Mahler volume due to Bourgain
and Milman.  We also generalize the bottleneck conjecture in the
context of indefinite orthogonal geometry and prove some special cases
of the generalization.}


\maketitlepage

Let $V$ be an $n$--dimensional vector space and let $V^*$ be the dual vector
space.  We denote the usual inner product between $V$ and $V^*$ by
$\bracket{\,\cdot\,,\,\cdot\,}$. If $K \subset V$ is a centrally symmetric convex body
centered at the origin, then there is a convex body
    $$K^\dual = \bigl\{\vy \in V^* \mid \bracket{K,\vy} \subseteq [-1,1]\bigr\}$$
called the {\em dual} or {\em polar} body of $K$.  The {\em Mahler volume} of $K$ is
defined as
    $$M(K) = \Vol K \times K^\dual = (\Vol K)(\Vol K^\dual).$$
Here $V$ and $V^*$ are given dual volume structures, or for the first
expression, the natural volume structure on $V \times V^*$ suffices.

The Mahler volume arises in the geometry of numbers and in functional
analysis.  By construction it is invariant under the action of $\GL(V)$ on $K$.
For fixed $V$, the space of symmetric convex bodies divided by the action of
$\GL(V)$ is compact in the Hausdorff topology, and $M(K)$ is continuous under
this action.  Consequently $M(K)$ has a finite maximum and a non-zero minimum
in each dimension. The maximum and minimum of $M(K)$ are interesting objects of
study in asymptotic convex geometry:

\begin{theorem}[Santal\'o] In a fixed vector space $V$, $M(K)$ is uniquely
maximized by ellipsoids. \label{th:santalo}
\end{theorem}

Let $C_n$ be the standard unit cube and let $B_n$ be the round unit ball, both
in $\R^n$. The polar body $C_n^\dual$ is the standard {\em cross polytope},
while  obviously $B_n^\dual = B_n$.

\begin{conjecture}[Mahler] For convex bodies $K$ in $n$ dimensions with $n$
fixed, the volume $M(K)$ is minimized by the cube $C_n$. \label{c:mahler}
\end{conjecture}

\Co{c:mahler} is considered harder than \Thm{th:santalo} because a
cube has much less symmetry than an ellipsoid. Moreover, $M(K)$ cannot be 
uniquely minimized when $K$ is a cube or a cross polytope, because there
are other polytopes with the same Mahler volume.  For example,
    $$M(C_{a+b}) = M(C_a \times C_b^\dual).$$
By contrast, \Thm{th:santalo} can be proved by an elegant symmetrization
argument \cite{Saint-Raymond:volume}.

Using methods from functional analysis, Bourgain and Milman \cite{BM:volume} proved an
asymptotic version of \Co{c:mahler}:

\begin{theorem}[Bourgain, Milman] There is a constant $c>0$ such that 
for any $n$ and any centrally-symmetric convex body $K$ of dimension $n$,
    $$M(K) \ge c^n M(B_n).$$ \label{th:bm}
\vspace{-\baselineskip}\end{theorem}

Although the proof technically constructs the constant $c$ (and although the
proof has been simplified \cite{Pisier:volume}), no good value for it is currently
known.  The author \cite{kuperberg:convex} proved the following:

\begin{theorem} If $K$ has dimension $n \ge 4$, then
    $$M(K) \ge (\log_2 n)^{-n} M(B_n).$$ \label{th:convex}
\vspace{-\baselineskip}\end{theorem}

\Thm{th:convex} has no arbitrary constants and therefore has some strength in
low dimensions, but it is obviously asymptotically weaker than
\Thm{th:bm}.

In this paper, we present a conjecture (\Co{c:bottleneck} below) 
which would produce a good value for the constant
$c$ in \Thm{th:bm}.  The conjecture also motivated the proof of
Theorem~\ref{th:convex}.

Let
\begin{align*}
K^+ & = \bigl\{(\vx,\vy) \in K \times K^\dual \mid \bracket{\vx,\vy} =  1\bigr\} \\
K^- & = \bigl\{(\vx,\vy) \in K \times K^\dual \mid \bracket{\vx,\vy} = -1\bigr\}
\end{align*}
and let $K^\diamond$ be the convex hull of $K^+ \cup K^-$. 

\begin{conjecture} For convex bodies $K$ in $n$ dimensions with $n$ fixed the
volume
    $$D(K) = \Vol K^\diamond$$
is uniquely minimized when $K$ is an ellipsoid. \label{c:bottleneck}
\end{conjecture}

\begin{fullfigure}{f:bottleneck}
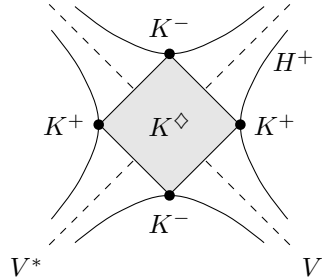
    {The geometry of $K^+$, $K^-$, and $K^\diamond$}
\pspicture(-2.1,-2.1)(2.1,1.7)
\psline[linestyle=dashed](-1.6,-1.6)(1.6, 1.6)
\psline[linestyle=dashed](-1.6, 1.6)(1.6,-1.6)
\rput{45}(0,0){
\pscurve( 2, .222)( 1.333, .333)( .666, .666)( .333, 1.333)( .222, 2)
\pscurve(-2, .222)(-1.333, .333)(-.666, .666)(-.333, 1.333)(-.222, 2)
\pscurve( 2,-.222)( 1.333,-.333)( .666,-.666)( .333,-1.333)( .222,-2)
\pscurve(-2,-.222)(-1.333,-.333)(-.666,-.666)(-.333,-1.333)(-.222,-2)
\psframe[fillstyle=solid,fillcolor=silver](-.666,-.666)(.666,.666)
\qdisk( .666, .666){2pt} \qdisk(-.666, .666){2pt}
\qdisk( .666,-.666){2pt} \qdisk(-.666,-.666){2pt}}
\rput(0,0){$K^\diamond$}
\rput(-1.4,0){$K^+$} \rput(1.4,0){$K^+$}
\rput(0,-1.3){$K^-$} \rput(0,1.3){$K^-$}
\rput(-1.9,-1.9){$V^*$} \rput(1.9,-1.9){$V$}
\rput(1.65,.85){$H^+$}
\endpspicture
\end{fullfigure}

We call \Co{c:bottleneck} the {\em bottleneck conjecture}, because the equation
$\bracket{\vx,\vy} = 1$ defines a hyperboloid sheet $H^+$ in $V \times V^*$
that resembles the flange of a bottle, while $K^+$ is a topological sphere in
$H^+$ that forms a neck.  \Fig{f:bottleneck} shows the geometry in the trivial
case $n=1$, which serves as a schematic for the higher-dimensional case. The
inclusion
    $$K^\diamond \subseteq K \times K^\dual$$
obviously implies the inequality
    $$D(K) \le M(K).$$

To see the strength of \Co{c:bottleneck}, consider these volume formulas:
\begin{align*}
\Vol C_n & = 2^n & \Vol C_n^\dual & = \frac{2^n}{n!} \\
\Vol B_n & = \frac{\pi^{n/2}}{(n/2)!} &
\frac{M(C_n)}{M(B_n)} & = \frac{(4/\pi)^n}{\binom{n}{n/2}}
\end{align*}
(Here $\frac{n}2! = \Gamma(\frac{n}2+1)$ when $n$ is half-integral.) The body
$B_n^\diamond$ is the convex hull of two orthogonal round $n$--balls of radius
$\sqrt{2}$ in $\R^{2n}$, so
    $$\Vol B_n^\diamond = (\Vol B_n)^2\frac{2^n}{\binom{2n}{n}}.$$
Consequently, if $\{K_n\}$ is any sequence of symmetric convex bodies with
$\dim K_n$\break$= n$, then \Co{c:mahler} implies that
    $$\lim_{n \to \infty} \sqrt[n]{\frac{M(K_n)}{M(B_n)}} \ge \frac2\pi$$
if the limit exists, while \Co{c:bottleneck} implies that
    $$\lim_{n \to \infty} \sqrt[n]{\frac{M(K_n)}{M(B_n)}} \ge \frac12$$
if the limit exists.

\section{Reformulations}
\label{s:reform}

The main purpose of this section is to introduce another conjecture which
implies \Co{c:bottleneck} and which may be equivalent.

\begin{conjecture} If $K \subset V$ is a centrally
symmetric convex body, then
    $$Q(\vecvol K^+),$$
the energy of the directed volume enclosed by $K^+$, is uniquely minimized when
$K$ is an ellipsoid. \label{c:energy}
\end{conjecture}

Here is an explanation of the terminology of \Co{c:energy}. The space $W = V
\times V^*$ has a symmetric bilinear form extending  the pairing of $V$ and
$V^*$ and such that
    $$\bracket{\vx,\vy} = 0$$
if $\vx$ and $\vy$ are both in $V$ or both in $V^*$.  (There is an even more
important antisymmetric, or symplectic, form that extends the pairing, but in
this article the symmetric extension is the relevant one.)  The function $Q$ is
the associated quadratic form on $W$ given by
    $$ Q(\vv) = \bracket{\vv,\vv}. $$
These forms have signature $(n,n)$, where $n$ is the dimension of $V$.  Both
the inner product and the quadratic form extend to the exterior algebra
$\Wedge^* W$ by the relation
    $$ Q(\omega_1 \wedge \omega_2) = Q(\omega_1)Q(\omega_2). $$
In this paper the quantity $Q(\omega)$ is called the {\em energy} of the tensor
$\omega$. The energy form $Q$ on the space $\Wedge^k W$ of  $k$--tensors has
signature
    $$ (\frac{a+b}2,\frac{a-b}2), $$
where
\begin{align*}
    a &= \binom{2n}{k} &
    b &= \begin{cases} (-1)^k \binom{n}{k/2} & \text {$k$ even} \\
        0 & \text{$k$ odd.} \end{cases}
\end{align*}

If $M \subset W$ is an oriented smooth $k$--manifold with boundary, it has a
{\em directed volume}
    $$\vecvol M \in \Wedge^k W.$$
If $M$ is the image of a smooth embedding
    $$\vf\co U \to W$$
of some domain $U \subset \R^k$, then the directed volume is given by an
integral formula:
    $$\vecvol M = \int_U d\vf = \int_U \frac{\d \vf}{\d x_1} \wedge
    \frac{\d \vf}{\d x_2} \wedge \ldots \wedge \frac{\d \vf}{\d x_k} d\vx.$$
By Stokes' theorem, $\vecvol M$ only depends on the boundary of $M$. If $N$ is
an oriented, closed $(k-1)$--manifold, we define the directed volume $\vecvol N$
{\em enclosed} by $N$ as the directed volume of any oriented $M$ with $\d M =
N$.

\subsection{\Co{c:energy} implies \Co{c:bottleneck}}
\label{s:implies}

The point of \Co{c:energy} is that the energy of the directed volume of $K$ is,
up to a constant factor, the volume of the region $K^\heart \subset K^\diamond$
enclosed by line segments that connect $K^+$ to $K^-$. The bodies $K^\heart$
and $K^\diamond$ could be identical for all $K$.  We will develop some
geometric properties of $K^-$ and $K^+$ to argue that
    $$ Q(\vecvol K^+) $$
is essentially an integral formula for the volume of $K^\heart$.

A vector $\vv \in W$ is {\em spacelike} if $Q(\vv) > 0$, {\em timelike} if
$Q(\vv) < 0$, and {\em null} if $Q(\vv) = 0$. A manifold in $W$ is {\em
spacelike} if all tangent vectors are spacelike; it is {\em timelike}
if all tangent vectors are timelike. There is a principle of transversality of
space and time: If $V^+$ is a spacelike vector subspace of $W$ and
$V^-$ is a timelike vector subspace, then
    $$ V^- \cap V^+ = \{\vec{0}\}. $$
Thus, any basis of $V^+$ and any basis of $V^-$ are linearly independent in
$W$.

Let $H^+$ and $H^-$ be the hypersurfaces defined by
    $$H^\pm = \bigl\{\vv \,\st\, Q(\vv) = \pm \frac12 \bigr\}.$$
Both hypersurfaces are diffeomorphic to $\R^n \times S^{n-1}$.   Pick some
ellipsoid $E \subset V$ centered at the origin. Then $E$ determines a
self-adjoint isomorphism
    $$\phi\co V \to V^*$$
such that
    $$ E = \bigl\{\vx \in V \bigm| \bracket{\vx,\phi(\vx)} \le 1\bigr\}. $$
Let $V^+$ and $V^-$ be the $n$--planes in $W$ defined by
    $$ V^\pm = \bigl\{(\vx,\pm \phi(\vx))\bigr\}. $$
Then
    $$ E^\pm = V^\pm \cap H^\pm. $$
The linear space $V^+$ is spacelike, while $V^-$ is timelike. The projection
of $H^+$ onto $V^+$ along $V^-$ consists of all points of $V^+$ except those
enclosed by $E^+$.  The composition of this linear projection with radial
projection onto $E^+$ is a convenient map
    $$ \pi^+\co H^+ \to E^+ $$
to $E^+$, which is a topological $(n-1)$--sphere.  Each fiber $\pi^{-1}(\vv)$ of
this map is a timelike section of $H^+$ which is isometric to hyperbolic
$n$--space.

As before, let $K$ be a symmetric convex body in $V$. For simplicity, assume
that both $K$ and $K^\dual$ are smooth.  For each point $\vx \in \d K$, there
is a unique $\vy \in \d K^\dual$, the outward normal of $\d K$ at $\vx$, such
that
    $$ \bracket{\vx,\vy} = 1. $$
Moreover, for each such $\vx$, the body $K$ has an osculating ellipsoid
$E(\vx)$, defined as the unique ellipsoid with the following three properties:
\begin{enumerate}
    \item $\vx$ lies in $\d E(\vx)$.
    \item $\vy$ is the outward normal of $E(\vx)$ at $\vx$.
    \item $\d E(\vx)$ has the same extrinsic curvature as $\d K$ at $\vx$.
\end{enumerate}
Equivalently, $E(\vx)^+$ and $K^+$ have the same tangent $(n-1)$--plane at the
point $(\vx,\vy)$.  The existence of $E(\vx)$ for each $\vx$ implies that $K^+$
is a spacelike manifold, \ie, that its tangent spaces are spacelike.  In fact,
for each $\vv \in K^+$, the $n$--plane spanned by $T_{\vv} K^+$ and $\vv$ is
spacelike.  Finally, the restriction of the projection $\pi^+$ to $K^+$ is a
homeomorphism between $K^+$ and $B_n^+$.

Let $J = K^+ * K^-$ be the topological join of $K^+$ and $K^-$. Explicitly,
    $$ J = \bigl(K^+ \times K^- \times [0,1]\bigr)/\!\sim $$
where the equivalence relation $\sim$ is given by
\begin{align*}
    (\vx,\vy_1,0) & \sim (\vx,\vy_2,0) &
    (\vx_1,\vy,1) & \sim (\vx_2,\vy,1).
\end{align*}
There is a natural map
    $$\vj\co J \to W$$
defined by
    $$\vj(\vx,\vy,t) = t\vx + (1-t)\vy.$$

In the following proposition and below, the adverb {\em almost} means  ``up to
a set of measure 0''.

\begin{proposition} The map $\vj$ is almost a smooth embedding. The set
$\vj(J)$ meets almost every ray from the origin in $W$ exactly once.
\label{p:almost}
\end{proposition}

\begin{proof} Let $S_W$ be the space of such rays, and let
    $$\pi_W\co J \to S_W$$
be the composition of $\vj$ with radial projection to $S_W$.  The space $J$ is
a smooth manifold except on $K^+$ and $K^-$, where it is merely a Lipschitz
manifold.  Let $\vx \in K^+$ and $\vy \in K^-$.  By the space-time
transversality principle, the vectors and tangent spaces $\vx$, $T_{\vx} K^+$,
$\vy$, and $T_{\vy} K^-$ are linearly independent.  Thus, the map $\pi$ has
positive Jacobian at each point $(\vx,\vy,t) \in J$ with $0 < t < 1$, because
the derivative matrix can be explicitly expressed in terms of $\vx$, $\vy$, and
bases for $T_{\vx} K^+$ and $T_{\vy} K^-$.  In other words, $\pi$ is a local
diffeomorphism away from $K^+$ and $K^-$.  The map $\pi$ is Lipschitz on $K^+$
and $K^-$ themselves, which implies that $\pi_W(K^+)$ and $\pi_W(K^-)$ are sets of
measure zero.

The degree of the map $\pi_W$ is both an integer and continuous as a function of
$K$.  It follows that the degree is 1, since that is its value when $K$ is an
ellipsoid.  Thus $\pi$ is almost a diffeomorphism, as desired.
\end{proof}

We conjecture that $\pi_W$ is a homeomorphism (without excepting a set
of measure zero).

As mentioned above, $K^\heart$ is defined as the region in $W$ enclosed
by $\vj(J)$.  By \Prop{p:almost}, $K^\heart$ is almost starlike.

Let $\vx \in K^+$ and let $P(\vx)$ be a tangent infinitesimal parallelepiped at
$\vx$.  Let $\vy \in K^-$ and define $P(\vy)$ likewise.  Let $P(\vx,\vy)$ be
the semi-infinitesimal polytope which is the convex hull of $P(\vx)$, $P(\vy)$,
and the origin. If the directed  volume of $P(\vx)$ is $d\vx$ and the directed
volume of $P(\vy)$ is $d\vy$, then the volume of $P(\vx,\vy)$ is
    $$\frac{1}{\binom{2n}{n}} \vx \wedge \vy \wedge d\vx \wedge d\vy.$$
The body $K^\heart$ is disjoint union of all $P(\vx,\vy)$ as $\vx$ and $\vy$
vary, and by Proposition~\ref{p:almost}, they are almost disjoint. 
Consequently
    $$\Vol K^\heart = \int_{K^+} \int_{K^-} \frac{1}{\binom{2n}{n}}
    \vx \wedge \vy \wedge d\vx \wedge d\vy.$$
This equation factors as
\begin{equation}
    \binom{2n}{n} \Vol K^\heart = \Bigl(\int_{K^+} \vx \wedge d\vx\Bigr)
    \wedge \Bigl(\int_{K^-} \vy \wedge d\vy\Bigr). \label{e:factors}
\end{equation}


Let $L^+$ be the union of line segments from $K^+$ to the origin and
let $L^-$ be the analogous cone over $K^-$.  Then
\begin{equation}
    \vecvol K^\pm = \vecvol L^\pm = \int_{K^\pm} \vx \wedge d\vx
    \label{e:integral}
\end{equation}
by decomposition into infinitesimal cones.  Thus, equation~(\ref{e:factors})
further simplifies to
\begin{align}
    \binom{2n}{n}\Vol K^\heart &= (\vecvol L^+) \wedge (\vecvol L^-) \notag \\
    & = (\vecvol K^+) \wedge (\vecvol K^-). \label{e:lk}
\end{align}
Finally, the linear map
    $$\sigma\co W \to W$$
defined by
    $$\sigma(\vx,\vy) = (-\vx,\vy)$$
for $\vx \in V$ and $\vy \in V^*$
sends $K^+$ to $K^-$ and negates the quadratic form $Q$. Both $\sigma$ and $Q$
extend to the exterior algebra $\Wedge^* W$. Functoriality of directed volume then implies that
\begin{equation}
    \vecvol K^- = \sigma \vecvol K^+. \label{e:sigmak}
\end{equation}

If $\ve_1,\ldots,\ve_n$ is a basis for $V$, and if $\phi$ is 
a self-adjoint isomorphism from $V$ to $V^*$ (as defined
previously), then
    $$\ve_1+\phi(\ve_1),\ve_2+\phi(\ve_2),\ldots,\ve_n+\phi(\ve_n)$$
is a basis for $V^+$ (also defined previously).  Then because
$\phi$ is self-adjoint, the wedge
product
    $$\omega = (\ve_1+\phi(\ve_1)) \wedge (\ve_2+\phi(\ve_2)) \wedge \ldots
    \wedge (\ve_n+\phi(\ve_n))$$
satisfies the identity
\begin{equation}
    \bracket{\omega,\nu} = \sigma(\omega) \wedge \nu \label{e:sigmawedge}
\end{equation}
for an arbitrary $n$--tensor $\nu$.  (It is easy to verify 
this identity with an explicit calculation in the
representative case where $V$ is $\R^n$ with the standard
basis and $\phi$ is the identity.)  Because of the system
of osculating ellipsoids for $K$, and because of equation
\eqref{e:integral}, $\vecvol K^+$ is a linear combination of
such tensors $\omega$, which means that it satisfies
equation \eqref{e:sigmawedge} as well.  In particular,
    $$\vecvol K^+ \wedge \sigma(\vecvol K^+) = Q(\vecvol K^+).$$
Combining this identity with equations
\eqref{e:lk} and \eqref{e:sigmak} yields
$$\binom{2n}{n}\Vol K^\heart = Q(\vecvol K^+).$$
Since $K^\heart$ is always contained in $K^\diamond$, and since they coincide
when $K$ is an ellipsoid, this final expression shows that \Co{c:energy} implies
\Co{c:bottleneck}, as desired.

\subsection{A generalization}

There is a plausible generalization of \Co{c:energy} to $a+b$ dimensions, by
which we mean a vector space $V$ with an inner product of signature $(a,b)$.
Let $Q$ be the associated quadratic form. Let
    $$H^+ = \bigl\{\vx \in V\, \st\, Q(\vx) = 1\bigr\}$$
be the positive unit hyperboloid sheet associated to $Q$.  (Note that $H^+$ is
now slightly different, because it was previously the level set
$Q^{-1}(1/2)$.)  Also for convenience endow $V$ with a volume form relative to
which the inner product has determinant $(-1)^b$.

\begin{conjecture} Let $H^+$ be the positive unit hyperboloid of a non-singular
quadratic form $Q$ on a vector space $V$ with signature $(a,b)$. Let
$N$ be a spacelike submanifold of $H^+$ whose inclusion into $H^+$ is a
homotopy equivalence. Then $Q(\vecvol N)$, the energy of the directed volume
enclosed by $N$, is uniquely minimized when $N$ is the intersection of $Q$ with
an $a$--plane in $V$ containing the origin. \label{c:general}
\end{conjecture}

Call a manifold $N$ as defined in \Co{c:general} a {\em neck}. \Co{c:energy} is
the special case of \Co{c:general} when $a=b$, and only for those necks which
can be realized as $K^+$ for some convex body $K$.

We could even more generally ask to minimize the inner product
    $$\bigl\langle\vecvol N_1,\vecvol N_2\bigr\rangle$$
for two different spacelike necks $N_1$ and $N_2$.  Or we could minimize the
wedge product
    $$\vecvol N^+ \wedge \vecvol N^-$$
for a spacelike neck $N^+$ in $H^+$ and a timelike neck $N^-$ in $H^-$. (The
wedge product can be interpreted as a number using the volume form on $V$.) In
the author's opinion, \Co{c:general} is a natural starting point for this
family of questions.

\section{Proofs in marginally indefinite cases}
\label{s:marginal}

In this section we will prove \Co{c:general} in the four least indefinite
cases:  $1+n$, $n+1$, $n+2$, and $2+n$ dimensions.  Note that in an
$(a+b)$--dimensional vector space $V$, the set of spacelike $a$--planes is
contractible, so we can consistently orient them.  Likewise we can consistently
orient timelike $b$--planes.  For convenience, we choose orientations which are
consistent with the orientation of $V$ induced by its volume form.

\subsection{Dimensions $1+n$ and $n+1$}

The first case, $1+n$ dimensions, is elementary. In this case $H^+$ is a
hyperboloid with two sheets and $N$ consists of a pair of points $\vx$ and
$\vy$, one on each sheet.  We can assume that $\vx$ is a positive vector and
$\vy$ is a negative vector. The directed volume of $N$ is then
    $$\vecvol N = \vx - \vy,$$
which is the sum of two positive unit spacelike vectors $\vx$ and $-\vy$. It is
elementary that the sum is shortest when they are parallel.  (Indeed, if we
switch space with time, this is the simplest case of the twin paradox in
special relativity.) This is equivalent to the condition that $N$ is centered
at the origin, the only thing to prove in this case.

\newgray{gray5}{.5}
\newgray{gray6}{.6}
\newgray{gray7}{.7}
\newgray{gray8}{.8}
\newgray{gray9}{.9}
\psset{hatchwidth=.4pt,hatchsep=3pt}
\begin{fullfigure}{f:hole}
    {$\pi_W(N)$ rings the hole of $\pi_W(H^+)$}
\pspicture(-3,-1.7)(3,1.7)
\pscircle[linestyle=none,fillstyle=crosshatch*,hatchcolor=gray9](0,0){1.7}
\pscircle[linestyle=none,fillstyle=crosshatch*,hatchcolor=gray8](0,0){1.5}
\pscircle[linestyle=none,fillstyle=crosshatch*,hatchcolor=gray7](0,0){1.3}
\pscircle[linestyle=none,fillstyle=crosshatch*,hatchcolor=gray6](0,0){1.1}
\pscircle[linestyle=none,fillstyle=crosshatch*,hatchcolor=gray5](0,0){.9}
\pscircle[fillstyle=solid,fillcolor=white](0,0){.7}

\psccurve(1.05;45)(1.35;90)(1.5;105)(1.45;155)(1.35;175)(1;235)(1;285)
    (1;315)(1;355)
\rput[l](1.6,.92){$\pi_W(H^+)$}
\rput[r](-2.3,.0){$\pi_W(N)$}
\psline{->}(-2.2,0)(-1.4,0)
\endpspicture
\end{fullfigure}

The second case, $n+1$ dimensions, is instructive for the last two cases, which
are more difficult.  Let $v_n$ be the volume of the unit ball in $\R^n$. Let
$W$ be a spacelike $n$--plane passing through the origin and let
    $$S = W \cap H^+$$
be the unit sphere in $W$.  Let
    $$\pi_W\co V \to W$$
be the orthogonal projection onto $W$, and let
    $$\pi_S\co H^+ \to S$$
be the radial projection onto $S$, generalizing the map $\pi^+$ of
\Sec{s:implies}.  By the argument of \Sec{s:implies}, $\pi_S$, if restricted to
$N$, is a homeomorphism. Equivalently, $\pi_W(N)$ is starlike.  At the same
time, $\pi_W(H^+)$ is the complement of $S$.  Consequently the area enclosed
by $\pi_W(N)$ is at least $v_n$, the volume enclosed by $B$, because
$\pi_W(N)$ must go around the hole in $\pi_W(H^+)$, as indicated
in \Fig{f:hole}.

Thus for any spacelike $n$--plane $W$, the component of $\vecvol N$ which is
orthogonal to $W$ is at least $v_n$.  This implies that $\vecvol N$ is dual to
a timelike vector.  If we choose an orthonormal basis
    $$\ve_1,\ve_2,\ldots,\ve_n$$
of $W$ and extend with a postive orthogonal unit timelike vector $\ve_{n+1}$,
$\vecvol N$ becomes the monomial tensor
    $$\vecvol N = c\, \ve_1 \wedge \ldots \wedge \ve_n.$$
Moreover, $c \ge v_n$, so by computation in this basis,
    $$Q(\vecvol N) \ge v_n^2.$$
The point is that in a suitable basis for $V$, the only
non-vanishing terms of $\vecvol N$ all have non-negative self inner product.

\subsection{Dimensions $n+2$ and $2+n$}

The third case, $n+2$ dimensions, requires a preliminary lemma about the
exterior square $\Wedge^2 V$ interpreted as a Lie algebra:
    $$\Wedge^2 V \cong \so(V) \cong \so(n,2).$$
Note that the first isomorphism is canonical, and that using this isomorphism,
    $$\bracket{X,Y} = -\frac12\Tr(XY).$$
Among the elements of $\so(V)$ there are spacelike and timelike rotations. 
Since the timelike planes are all oriented, the timelike rotations can be
divided into positive and negative.  Also say that an element of $\so(V)$ is
{\em elliptic} if it is a product of commuting spacelike and timelike rotations
(positive or negative).

\begin{lemma}[Paneitz] A convex combination of positive timelike rotations is 
elliptic. \label{l:paneitz}
\end{lemma}

Here are some comments about the results and terminology of Paneitz
\cite{Paneitz:causality,Paneitz:determination}. Among all convex cones in
$\so(V)$ which are invariant under conjugation, there is a unique minimal
closed cone $C_0$ and  a unique maximal cone $C_1$ (necessarily closed). 
Define the infinitesimal angle $d > 0$ of a rotation $R$ (either spacelike or
timelike) by the relation
    $$\Tr(R^2) = 2d^2.$$
Then according to Paneitz \cite[page 340]{Paneitz:causality}, the elements of
$C^{\text{int}}_0$ are precisely those that are a commuting product of a
positive timelike rotation by an angle $d_0$ and spacelike rotations by angles
$d_1,\ldots,d_k$ (necessarily $2k \le n$) such that
    $$d_0 > d_1+d_2+\ldots+d_k.$$
Every timelike rotation is of this form (with $k=0$), hence any convex
combination is as well.

Recall that an alternating $k$--tensor is {\em simple} if it is a wedge product
of vectors. For a general quadratic form $Q$ on $V$ of signature $(a,b)$, say
that a simple $k$--tensor in $\Wedge^k V$ is {\em spacelike} (respectively {\em
timelike}) if it is the wedge product of vectors that span a spacelike
$k$--plane (resp. a timelike $k$--plane). A spacelike simple $a$--tensor (resp. a
timelike simple $b$--tensor) is {\em positive} if its factors are positively
ordered relative to the orientation of the $a$--plane (resp. the $b$--plane) they
span.  Recall that the Hodge star operator on $k$--tensors is defined as the
unique linear operator
    $$*\co \Wedge^k V \to \Wedge^{n+2-k} V$$
such that
    $$*(\ve_1 \wedge \ve_2 \wedge \ldots \wedge \ve_k)
    = \ve_{k+1} \wedge \ve_{k+2} \wedge \ldots \wedge \ve_{n+2}$$
for any positively oriented orthonormal frame
    $$\ve_1,\ve_2,\ldots,\ve_{n+2}.$$
We will need two facts about the Hodge star operator: first, that 
    $$Q(*\omega) = (-1)^b Q(\omega)$$
for any tensor $\omega$, and second that $\omega$ is a positive, spacelike,
simple $a$--tensor if and only if $*\omega$ is a positive, timelike, simple
$b$--tensor.

In terms of $2$--tensors, \Lem{l:paneitz} says that a convex combination of
positive, timelike, simple 2--tensors can be expressed in the form
    $$d_0 \ve_0 \wedge \ve_1 + d_1 \ve_2 \wedge \ve_3 + \ldots +
    d_k \ve_{2k} \wedge \ve_{2k+1},$$
where the vectors $\ve_0,\ve_1,\ldots,\ve_{2k+1}$ are orthonormal, and $\ve_0$
and $\ve_1$ are timelike. In addition if the pair $(\ve_0,\ve_1)$ forms a
positive basis of the plane it spans, then $d_0$ is positive. We will need the
dual statement that a convex combination of positive, spacelike, simple
$n$--tensors can be expressed in the form
\begin{align}
    d_0 *\!(\ve_0 \wedge \ve_1) & + d_1 *\!(\ve_2 \wedge \ve_3)
    + d_2 *\! (\ve_4 \wedge \ve_5) + \ldots \notag \\
    & + d_k *\!(\ve_{2k} \wedge \ve_{2k+1}). \label{e:stdform}
\end{align}

Finally, if $N$ is a neck, then $\vecvol N$ is realized as a convex combination
of positive, spacelike, simple $n$--tensors by the obvious generalization of
equation~\eqref{e:integral}. Consequently $*\vecvol N$ can be expressed in the
form of expression \eqref{e:stdform}.  If $W$ is a spacelike $n$--plane spanned
by the vectors $\ve_2,\ldots,\ve_{2k+1}$, then the projection of $N$ encloses a
volume of at least $v_n$ by the idea illustrated in Figure~\ref{f:hole}. Thus
$$d_0 \ge v_n,$$
and 
$$Q(\vecvol N) = \sum_{i=0}^k d_i^2 \ge d_0^2 \ge v_n^2,$$
as desired.

\Co{c:general} is argued the same way in $2+n$ dimensions as in
$n+2$ dimensions,
except without the complication of applying Hodge duality.

\subsection{Trivial cases and open cases}

The case of $n+0$ dimensions is trivially true, since there is only one
candidate for the neck $N$.  The case of $0+n$ dimensions is vacuous.

The basic reason that the above arguments do not work in $a+b$ dimensions when
both $a$ and $b$ are at least 3 is that the space of  alternating $a$--tensors
is bigger than the Lie group $\SO(a,b)$. Asymptotically
    $$\dim \Wedge^a \R^{a+b}$$
grows exponentially in $\min(a,b)$, while
    $$\dim \SO(a,b)$$
grows quadratically.  The general $a$--tensor does not admit an orthonormal
basis such that all terms have positive energy.

\section{Local stability}
\label{s:local}

In this section we argue that a flat neck is a local minimum of the energy
$Q(\vecvol N)$ relative to the $C^1$ topology in $a+b$ dimensions.

Consider $\R^{a+b}$ together with the standard quadratic form $Q$ of signature
$(a,b)$ given by
    $$Q(\vx,\vy) = \vx \cdot \vx - \vy \cdot \vy,$$
using the standard dot products on $\R^a$ and $\R^b$. Let
$\bracket{\,\cdot\,,\,\cdot\,}$ be the associated bilinear form.  Let $S^{a-1}$ be the
standard unit $(a-1)$--sphere in the standard timelike $\R^a \subset \R^{a+b}$. 
The hyperboloid sheet $H^+$ is perpendicular to $\R^a$ at the sphere $S^{a-1}$. 
Given a $C^1$ function
    $$\vf\co \R^a \to \R^b,$$
let $N$ be the set 
    $$N = \Bigl\{\bigl(\vx\sqrt{1+f^2(\vx)},\vf(\vx)\bigr)\,\st\,
    \vx \in S^{a-1}\Bigr\}.$$
For suitable $\vf$, $N$ is a neck, and every neck $N$ can be  uniquely
expressed in this form.

Let $\ve_1,\ldots,\ve_{a+b}$ be the standard basis of $\R^{a+b}$.
Given a linear map
    $$L\co \R^a \to \R^b,$$
we define an alternating $a$--tensor
$$\Psi(L) = \sum_{k=1}^a (-1)^{k+1} L(\ve_k) \wedge
    \ve_1 \wedge \ldots \wedge \widehat{\ve_k} \wedge \ldots \wedge \ve_a.$$
In other words, $\Psi$ is the natural linear transformation
$$\Psi\co \Hom(\R^a,\R^b) \to \Wedge^{a-1} \R^a \tensor \Wedge^1 \R^b
    \subset \Wedge^a \R^{a+b}$$
induced by the standard Hodge star operator on $\R^a$ and the standard dot
product on $\R^b$. Using this notation, if $\vf$ and its derivative $D\vf$ are
of order $\epsilon$, then
\begin{align}
    \vecvol N =\ & \Bigl(v_{a-1} + \int af^2(\vx)\,d\vx
        + o(\epsilon^2)\Bigr) \ve_1 \wedge \ldots \wedge \ve_a \notag \\
        &+ \int \Psi(D\vf(\vx))\,d\vx + o(\epsilon). \label{e:mess}
\end{align}
Here all integrals are over the sphere $S^{a-1}$, as before $v_{a-1}$ is the
volume enclosed by $S^{a-1}$, and the last term $o(\epsilon)$ consists of
monomials with at least two wedge factors $\ve_k$ with $k>a$.   If we set
    $$Q[\vf] = Q(\vecvol N),$$
then the first variational derivative of $Q$ at $\vf = 0$ vanishes by symmetry,
while the second variational derivative is given by
\begin{equation} \begin{array}{rccc}\displaystyle
    \left.\frac{\delta^2 Q}{ (\delta \vf)^2}\right|_0 = &
    \displaystyle \int af^2(\vx)\,d\vx\ & - &
    \displaystyle \left(\int \Psi(D\vf(\vx))\ d\vx\right)^2 \\
    \stackrel{\text{def}}{=} & A[\vf] & - &\ B[\vf] 
\end{array} \label{e:diff} \end{equation}
from equation~\eqref{e:mess}.  In the second line of equation~\eqref{e:diff},
we define the functional $A[\vf]$ to be the first term of the
first line and the functional $B[\vf]$ to be the second term.

We claim that the second variational derivative of $Q$ (equation
\eqref{e:diff}) is positive definite except for null directions given by the
action of the symmetry group $\SO(a,b)$.  These null directions correspond to
the variations $\vf$ which are linear. The general $\vf$ has
a harmonic expansion
    $$\vf = \vf_0 + \vf_1 + \vf_2 + \ldots,$$
where $\vf_k$ is given by a degree $k$ polynomial which is orthogonal to
lower-degree polynomials on the sphere $S^{a-1}$.  The functional
$A$ is proportional to the $L^2$ norm of $\vf$:
\begin{equation}
    A[\vf] = a||\vf||^2 = a||\vf_0||^2 + a||\vf_1||^2 +
    a||\vf_2||^2 + \ldots. \label{e:a}
\end{equation}
On the other hand, the functional $B$ is a quadratic function composed
with the linear transformation
$$\vf \mapsto \int \Psi(D\vf(\vx))\ d\vx.$$
This transformation is equivariant under $\SO(a) \times \SO(b)$, the stabilizer
in $\SO(a,b)$ of the flat neck $S^{a-1}$.  Its target is the irreducible
representation $\Wedge^{a-1} \R^a \tensor \Wedge^1 \R^b$.  Therefore it must
annihilate all terms of the harmonic expansion of $\vf$ except for $\vf_1$, the
sole term which lies in an isomorphic summand of the $L^2$ completion of the
function space $C^1(S^{a-1},\R^b)$. In other words,
\begin{equation}
    B[\vf] = B[\vf_1] = c||\vf_1||^2 \label{e:b}
\end{equation}
for some constant $c$.  This constant $c$ can be determined by noting
that if $\vf$ is linear, \ie, $\vf = \vf_1$, then 
    $$A[\vf] - B[\vf] = 0,$$
because then $\vf$ represents an infinitesimal motion of the neck given by the
action of the Lie algebra $\so(a,b)$.  Consequently $c=a$. Subtracting
equation~\eqref{e:b} from equation~\eqref{e:a}, we obtain
    $$\left.\frac{\delta^2 Q}{ (\delta \vf)^2}\right|_0 =
    a||\vf_0||^2 + a||\vf_2||^2 + a||\vf_3||^2 + a||\vf_4||^2 + \ldots.$$
Thus the second variational derivative has the desired positivity
property.

\section{Acknowledgements}

The author would like to thank Krystyna Kuperberg and W{\l}odzimierz Kuperberg
for corrections as well as Bruce Kleiner and Frank Morgan for their interest in
the work.

\end{document}